\documentclass[11pt]{amsart}
\usepackage[T1]{fontenc}
\usepackage{lmodern,amssymb}
\usepackage[expansion=false]{microtype}
\usepackage{xurl}
\usepackage[hidelinks]{hyperref}

\newtheorem{theorem}{Theorem}[section]
\newtheorem{proposition}[theorem]{Proposition}
\newtheorem{corollary}[theorem]{Corollary}
\theoremstyle{remark}
\newtheorem{remark}[theorem]{Remark}
\newcommand{\Pow}{\mathcal P}
\newcommand{\ltL}{<_{\mathrm L}}
\newcommand{\leNE}{\leq_{\mathrm{NE}}}
\newcommand{\ZF}{\mathsf{ZF}}
\newcommand{\DC}{\mathsf{DC}}
\newcommand{\doi}[1]{\href{https://doi.org/#1}{doi:\nolinkurl{#1}}}

\title[The strict Lipschitz order on filters]
  {Well-foundedness of the strict Lipschitz order on countably complete filters}
\author{Gabriel Fernandes}
\address{Instituto de Ci\^encias Matem\'aticas e de Computa\c c\~ao,
  Universidade de S\~ao Paulo, S\~ao Carlos, SP 13566-590, Brazil}
\email{fernandes@icmc.usp.br}
\thanks{Gabriel Fernandes was supported in part by grant \#2025/09425-1,
  S\~ao Paulo Research Foundation (FAPESP)}
\author{Matheus Duzi}
\address{Instituto de Ci\^encias Matem\'aticas e de Computa\c c\~ao,
  Universidade de S\~ao Paulo}
\email{matheus.duzi.costa@usp.br}
\hypersetup{pdftitle={Well-foundedness of the strict Lipschitz order on countably complete filters},
  pdfauthor={Gabriel Fernandes and Matheus Duzi},
  pdfkeywords={countably complete filters, ultrafilters, Lipschitz order, Ketonen order}}

\begin{document}

\begin{abstract}
In $\ZF+\DC$, the strict Lipschitz order on proper countably complete
filters on any fixed ordinal is well-founded. The proof uses transfinite
recursion and closure under countable intersections to derive a
contradiction from a descending chain. Restricting each filter to the
least ordinal it contains yields a set-like well-founded order that
compares filters on different ordinals. In $\ZF$, two ultrafilters on
the same ordinal are equal if each is nonexpansively reducible to the
other. These results answer both parts of Question~9.2.11 in
\emph{The Ultrapower Axiom}.
\end{abstract}

\maketitle

\section{Introduction}

The relations studied here bring together two approaches to comparing
complexity. In descriptive set theory, Wadge reducibility compares
sets by continuous inverse images: an equality $A=f^{-1}[B]$ transforms
the question of membership in $A$ into that of membership in $B$.
This gives a finer analysis than membership in the usual definability
classes. Lipschitz reducibility restricts the available transformations
to nonexpansive maps (see \cite[Introduction and \S I.A]{Wadge}).

The associated games explain this restriction. Continuous reductions
allow the responding player to wait for more input before producing
the next output, whereas nonexpansive reductions require an immediate
response at each round. A super-Lipschitz reduction must produce the
current output using only earlier input. On $2^\omega$ this is the
contraction condition for the usual ultrametric. These restrictions on
coordinate dependence also make sense on $2^\delta$ for an arbitrary
ordinal, independently of a real-valued metric.

The ultrafilter orders have a complementary origin in the analysis of
large cardinals and canonical inner models. Mitchell's work on coherent
sequences of normal measures extended the analysis of $L[U]$ to models
with prescribed collections of measures \cite{Mitchell}. The Mitchell
order compares normal ultrafilters by membership in ultrapowers, and
its ordinal ranks measure their complexity \cite[\S1]{Goldberg}.
Ketonen's original comparison arose in a proof of his characterization
of strong compactness using complete uniform ultrafilters on regular
cardinals \cite{Ketonen}. In that proof, minimal weakly normal
ultrafilters provide witnesses to strong compactness (see
\cite[\S3.1.1]{GoldbergBook} for the historical account).

Goldberg's generalization of the Lipschitz order connects these
viewpoints: an ultrafilter on $\delta$ is itself a subset of
$\Pow(\delta)$, so its complexity can be compared by reductions.
Moreover, the reduction games become games of ordinal length.
This provides a link between ultrapower comparison and determinacy
\cite[\S5.2]{Goldberg}.
For arbitrary subsets, a single contraction may reduce a set to itself.
The strict relation therefore requires reductions both to the target
and to its complement, so that the fixed point argument excludes
self-comparison \cite[\S5.1]{Goldberg}.

The strict Lipschitz order is defined on subsets of $2^\delta$ for any
ordinal $\delta$ \cite[\S5]{Goldberg}. On countably complete ultrafilters
it extends the strict Ketonen order, and the two coincide under the
Ultrapower Axiom \cite[Theorem~5.13 and Lemma~5.16]{Goldberg}. The strict
Lipschitz order agrees with the Mitchell order on normal ultrafilters
and is therefore well-founded there \cite[Proposition~5.21]{Goldberg}.
It is also well-founded on Dodd sound ultrafilters
\cite[Corollary~4.3.28]{GoldbergBook}.

The Ketonen and Lipschitz orders can differ even under the Weak
Ultrapower Axiom \cite{Kaplan}. The discussion following Definition~2.3
of that paper also raises the well-foundedness question.

We address both parts of Question~9.2.11 in \cite{GoldbergBook}.
\begin{quote}
Must the Lipschitz order on countably complete ultrafilters be
wellfounded? Can there be Lipschitz equivalent (countably complete)
ultrafilters that are not equal?
\end{quote}

Theorem~\ref{thm:main} answers the first question affirmatively in
$\ZF+\DC$ for proper countably complete filters. Corollary~\ref{cor:nonexpansive}
answers the second negatively in $\ZF$ for arbitrary ultrafilters on
the same ordinal.

Section~\ref{sec:global-order} extends the comparison to filters on
different ordinals and proves that the resulting order is set-like and
well-founded on canonical presentations of countably complete filters.
Section~\ref{sec:incomplete-ultrafilters} shows that every countably
incomplete ultrafilter starts an infinite descending sequence on the
same ordinal.

\section{The strict order}

We identify $2^\delta$ with $\Pow(\delta)$. A map
$f:\Pow(\delta)\to\Pow(\delta)$ is \emph{super-Lipschitz} if for every
$\alpha<\delta$ and $X,Y\subseteq\delta$,
\[
 X\cap\alpha=Y\cap\alpha\quad\Longrightarrow\quad
 (\alpha\in f(X)\iff\alpha\in f(Y)).
\]
Such maps are called contractions in \cite[Definition~5.8]{Goldberg}.
Equivalently, for each $\alpha<\delta$ the set $f(X)\cap(\alpha+1)$
depends only on $X\cap\alpha$
\cite[discussion preceding Definition~2.3]{Kaplan}.

For $D,E\subseteq\Pow(\delta)$, write $D\ltL E$ if there are
super-Lipschitz maps $p,q:\Pow(\delta)\to\Pow(\delta)$ such that for
every $X\subseteq\delta$,
\begin{equation*}
 p(X)\in E\iff X\in D,
 \qquad
 q(X)\in E\iff X\notin D.
\end{equation*}
The map $p$ is a positive reduction and $q$ is a negative reduction.
Equivalently, $D$ reduces to both $E$ and $\Pow(\delta)\setminus E$
\cite[Definition~5.10]{Goldberg}. If $E$ is an ultrafilter, a positive
reduction $p$ supplies a negative reduction via $q(X)=\delta\setminus p(X)$.
Thus on ultrafilters the positive reduction alone suffices, as in
\cite[Definition~2.3]{Kaplan}.

Super-Lipschitz maps are closed under composition because they preserve
agreement below each $\alpha$. This gives transitivity of $\ltL$:
if $p,q$ witness $D\ltL E$ and $r$ is a positive witness for $E\ltL H$,
then $r\circ p$ and $r\circ q$ witness $D\ltL H$. Irreflexivity follows
from the fixed point property \cite[Theorem~5.9 and Lemma~5.11]{Goldberg}:
given a super-Lipschitz map $q$, recursively include $\alpha$ in $X$
exactly when $\alpha\in q(X\cap\alpha)$. Then $q(X)=X$. A negative
reduction for $D\ltL D$ would therefore give $X\in D\iff X\notin D$.

\begin{theorem}[$\ZF+\DC$]\label{thm:main}
For every ordinal $\delta$, the strict Lipschitz order on proper
countably complete filters on $\delta$ is well-founded. In particular,
it is well-founded on countably complete ultrafilters on $\delta$.
\end{theorem}

\begin{proof}
Suppose a nonempty set of proper countably complete filters on $\delta$
has no $\ltL$-minimal element. By dependent choice there is a sequence
$(D_n)_{n<\omega}$ in this set with $D_{n+1}\ltL D_n$ for each $n$.
Countable choice, which follows from dependent choice, gives positive
and negative reductions $p_n,q_n$ for these comparisons. For $n<m$ put
\[
 g_{n,m}=p_n\circ\cdots\circ p_{m-2}\circ q_{m-1},
\]
with no positive factors when $m=n+1$.
Every $g_{n,m}$ is super-Lipschitz and satisfies
\begin{equation}\label{eq:negative}
 g_{n,m}(X)\in D_n\quad\Longleftrightarrow\quad X\notin D_m.
\end{equation}

Define sets $B_\alpha\subseteq\omega$ by recursion on $\alpha<\delta$.
At stage $\alpha$, let $A_m^\alpha=\{\beta<\alpha:m\in B_\beta\}$ for
each $m<\omega$, and put
\[
 B_\alpha=\left\{n<\omega:\alpha\in\bigcap_{m>n}g_{n,m}(A_m^\alpha)\right\}.
\]
For each $m$ let $A_m=\{\alpha<\delta:m\in B_\alpha\}$. Fix $n<\omega$
and $\alpha<\delta$. For each $m>n$, we have $A_m\cap\alpha=A_m^\alpha$.
Since $g_{n,m}$ is super-Lipschitz, $\alpha$ belongs to
$g_{n,m}(A_m^\alpha)$ exactly when it belongs to $g_{n,m}(A_m)$.
It follows that
\[
 \begin{aligned}
  \alpha\in A_n
   &\iff n\in B_\alpha\\
   &\iff \alpha\in\bigcap_{m>n}g_{n,m}(A_m^\alpha)\\
   &\iff \alpha\in\bigcap_{m>n}g_{n,m}(A_m).
 \end{aligned}
\]
As this holds for every $\alpha<\delta$, we obtain
\begin{equation}\label{eq:intersection}
 A_n=\bigcap_{m>n}g_{n,m}(A_m)\qquad(n<\omega).
\end{equation}

By countable completeness and upward closure of $D_n$, the intersection
in \eqref{eq:intersection} belongs to $D_n$ exactly when each of its
terms does. Using \eqref{eq:negative}, we obtain
\begin{equation}\label{eq:tail}
 A_n\in D_n\quad\Longleftrightarrow\quad
 (\forall m>n)\ A_m\notin D_m.
\end{equation}
If $A_n\in D_n$ for some $n$, then $A_m\notin D_m$ for every $m>n$.
In particular, $A_{n+1}\notin D_{n+1}$, while \eqref{eq:tail} at $n+1$
gives $A_{n+1}\in D_{n+1}$. Hence $A_n\notin D_n$ for every $n$.
Applying \eqref{eq:tail} at $0$ now gives $A_0\in D_0$, a contradiction.
\end{proof}

\begin{remark}
The truth conditions in the paradox described in \cite{Yablo} have
the same form as the equivalences in \eqref{eq:tail}. In the recursive
step, the definition at each ordinal $\alpha<\delta$ depends only on
what has been defined at ordinals below $\alpha$. The same construction
works for filters on a set $I$ equipped with a well-founded strict
partial order $<$: at each $i\in I$, the definition depends only on
elements $j\in I$ with $j<i$. The corresponding maps
$f:\Pow(I)\to\Pow(I)$ must satisfy
\[
 i\in f(X)\quad\Longleftrightarrow\quad
 i\in f(X\cap\{j\in I:j<i\})
 \qquad(i\in I,\ X\subseteq I).
\]
\end{remark}

A map $f:\Pow(\delta)\to\Pow(\delta)$ is \emph{nonexpansive} if
for every $\alpha<\delta$ and all $X,Y\subseteq\delta$,
\[
 X\cap(\alpha+1)=Y\cap(\alpha+1)\quad\Longrightarrow\quad
 f(X)\cap(\alpha+1)=f(Y)\cap(\alpha+1).
\]
Equivalently, for each $\alpha<\delta$ and $X\subseteq\delta$,
we have $\alpha\in f(X)$ if and only if
$\alpha\in f(X\cap(\alpha+1))$.
This is the Lipschitz condition of \cite[Definition~3.4.24]{GoldbergBook}.
For $\delta=\omega$ it is the $1$-Lipschitz condition for the usual
ultrametric on $2^\omega$.
Write $U\leNE W$ if there is a nonexpansive map $f$ such that
$f^{-1}[W]=U$.

\begin{corollary}[$\ZF$]\label{cor:nonexpansive}
For any ultrafilters $U$ and $W$ on the same ordinal,
\[
 U\leNE W \quad\Longleftrightarrow\quad U=W\text{ or }U\ltL W.
\]
Thus nonexpansive reducibility is antisymmetric, its equivalence
classes are singletons, and its strict part is $\ltL$.
\end{corollary}

\begin{proof}
Suppose $U\ne W$ and $f$ witnesses $U\leNE W$. Choose $S\in U\setminus W$
and put $T=\delta\setminus S\in W$. Define
\[
 g(X)=f(X\cap S)\cap T\qquad(X\subseteq\delta).
\]
Since $S\in U$ and $T\in W$, closure under finite intersections and
upward closure give, for every $X\subseteq\delta$,
\[
 g(X)\in W\iff f(X\cap S)\in W\iff X\cap S\in U\iff X\in U.
\]
Thus $g$ still reduces $U$ to $W$.

Fix $\alpha<\delta$ and suppose $X\cap\alpha=Y\cap\alpha$. If $\alpha\in S$,
then $\alpha\notin T$, so both $g(X)$ and $g(Y)$ omit $\alpha$.
If $\alpha\in T$, intersecting with $S$ gives
\[
 (X\cap S)\cap(\alpha+1)
   =(X\cap\alpha)\cap S
   =(Y\cap\alpha)\cap S
   =(Y\cap S)\cap(\alpha+1).
\]
By nonexpansiveness, $f(X\cap S)$ and $f(Y\cap S)$ agree at $\alpha$,
and intersecting with $T$ preserves that agreement. In both cases
$g(X)$ and $g(Y)$ agree at $\alpha$, so $g$ is super-Lipschitz.
Since $W$ is an ultrafilter, this gives $U\ltL W$.

Conversely, the identity witnesses $U\leNE U$, and every super-Lipschitz
reduction is nonexpansive. If $U\ne W$ and each ultrafilter reduces
nonexpansively to the other, transitivity of $\ltL$ yields $U\ltL U$,
a contradiction.
\end{proof}

\section{Comparing filters on different ordinals}\label{sec:global-order}

To keep the construction self-contained, we first prove the basic
properties of extension, restriction, and strict Lipschitz comparison
for filters on ordinals.
All filters in this section are proper. For a filter $D$ on an ordinal
$\eta$ and an ordinal $\delta\geq\eta$, define its \emph{extension to
$\delta$} by
\[
 D^\delta=\{X\subseteq\delta:X\cap\eta\in D\}.
\]

If $\sigma\leq\eta$ and $\sigma\in D$, define the \emph{restriction to
$\sigma$} by
\[
 D\upharpoonright\sigma=D\cap\Pow(\sigma).
\]

The following are elementary properties of filter projection and change of
space (see \cite[Definitions~3.2.1--3.2.2]{GoldbergBook}).

\begin{proposition}\label{prop:extension-restriction}
Let $D$ be a proper filter on $\eta$, and let $\sigma\leq\eta\leq\delta
\leq\gamma$.
\begin{enumerate}
\item $D^\delta$ is a proper filter on $\delta$. Extension preserves
countable completeness and the ultrafilter property, and
$D^\eta=D$ and $(D^\delta)^\gamma=D^\gamma$.
\item If $\sigma\in D$, then $D\upharpoonright\sigma$ is a proper
filter on $\sigma$ and inherits countable completeness and the
ultrafilter property from $D$. Moreover,
\[
 X\in D\quad\Longleftrightarrow\quad
 X\cap\sigma\in D\upharpoonright\sigma
 \qquad(X\subseteq\eta).
\]
Consequently, $D=(D\upharpoonright\sigma)^\eta$.
\end{enumerate}
\end{proposition}

For a filter $D$ on $\eta$, let $\rho(D)$ be the least ordinal in $D$
and define its \emph{canonical presentation} by
\[
 \rho(D)=\min\{\sigma\leq\eta:\sigma\in D\},
 \qquad
 \operatorname{can}(D)=
 \bigl(\rho(D),D\upharpoonright\rho(D)\bigr)
\]
(for ultrafilters, these canonical representatives appear in
	\cite[Definition~3.2.5 and Lemmas~3.2.6--3.2.7]{GoldbergBook}).
The minimum exists because $\eta\in D$, and $\rho(D)>0$ because
$\varnothing\notin D$. More generally, a pair $(\rho,D)$ is
\emph{canonical} if $D$ is a proper filter on $\rho$ and
$\xi\notin D$ for every $\xi<\rho$. By
Proposition~\ref{prop:extension-restriction}, passing to a canonical
presentation preserves countable completeness and the ultrafilter property.
For $\delta\geq\eta$, the definition of extension gives
$\rho(D^\delta)=\rho(D)$ and
$D^\delta\upharpoonright\rho(D)=D\upharpoonright\rho(D)$.
Thus extension leaves the canonical presentation unchanged.

For filters $D,E$ on the same ordinal $\eta$, write $D\ltL^\eta E$
for the relation $D\ltL E$ of Section~2.
For canonical pairs $(\rho,D)$ and $(\sigma,E)$, put
$\tau=\max\{\rho,\sigma\}$ and define the \emph{global strict Lipschitz
order} by
\begin{equation}\label{eq:global-L}
 (\rho,D)\ltL^{\mathrm{glob}}(\sigma,E)
 \quad\Longleftrightarrow\quad
 D^\tau\ltL^\tau E^\tau.
\end{equation}
Explicitly, this means that there are super-Lipschitz maps
$p,q:\Pow(\tau)\to\Pow(\tau)$ such that, for every $X\subseteq\tau$,
\[
 \begin{aligned}
  p(X)\cap\sigma\in E&\quad\Longleftrightarrow\quad X\cap\rho\in D,\\
  q(X)\cap\sigma\in E&\quad\Longleftrightarrow\quad X\cap\rho\notin D.
 \end{aligned}
\]
To compare filters on arbitrary ordinals, apply $\ltL^{\mathrm{glob}}$
to their canonical presentations. The use of a common ordinal is
analogous to the definition of the global Ketonen order
\cite[Definition~3.3.13]{GoldbergBook}.

\begin{proposition}\label{prop:order-extension}
Let $\beta\leq\delta$ and let $D,E$ be filters on $\beta$. Then
\begin{equation}\label{eq:order-extension}
 D\ltL^\beta E
 \quad\Longleftrightarrow\quad
 D^\delta\ltL^\delta E^\delta.
\end{equation}
\end{proposition}

\begin{proof}
For the forward implication, let $p,q$ witness $D\ltL^\beta E$ and
define
\[
 \widehat p(X)=p(X\cap\beta),\qquad
 \widehat q(X)=q(X\cap\beta)
 \qquad(X\subseteq\delta),
\]
regarding $p(X\cap\beta)$ and $q(X\cap\beta)$ as subsets of $\delta$.
Fix $r\in\{p,q\}$, $\alpha<\delta$, and $X,Y\subseteq\delta$ with
$X\cap\alpha=Y\cap\alpha$. If $\alpha<\beta$, then
$(X\cap\beta)\cap\alpha=(Y\cap\beta)\cap\alpha$, so
super-Lipschitzness of $r$ gives
$\alpha\in\widehat r(X)\iff\alpha\in\widehat r(Y)$.
If $\beta\leq\alpha$, then $\widehat r(X),\widehat r(Y)\subseteq\beta$,
so both sets omit $\alpha$. Thus $\widehat p,\widehat q$ are
super-Lipschitz. Moreover,
\[
 \begin{aligned}
  \widehat p(X)\in E^\delta
    &\iff X\cap\beta\in D\iff X\in D^\delta,\\
  \widehat q(X)\in E^\delta
    &\iff X\cap\beta\notin D\iff X\notin D^\delta.
 \end{aligned}
\]
Conversely, let $P,Q$ witness $D^\delta\ltL^\delta E^\delta$ and put
\[
 p(A)=P(A)\cap\beta,\qquad q(A)=Q(A)\cap\beta
 \qquad(A\subseteq\beta).
\]
Here $A$ is also regarded as a subset of $\delta$. Suppose
$A,B\subseteq\beta$ and $A\cap\alpha=B\cap\alpha$ for $\alpha<\beta$.
For each $R\in\{P,Q\}$, super-Lipschitzness gives
$\alpha\in R(A)\cap\beta\iff\alpha\in R(B)\cap\beta$.
Thus $p,q$ are super-Lipschitz. Their reduction properties follow from
\[
 \begin{aligned}
  p(A)\in E&\iff P(A)\in E^\delta\iff A\in D,\\
  q(A)\in E&\iff Q(A)\in E^\delta\iff A\notin D.
 \end{aligned}
\]
\end{proof}

Applying Proposition~\ref{prop:order-extension} to $D^\tau,E^\tau$
and using the composition law for extensions gives, for every
$\delta\geq\tau$,
\[
 (\rho,D)\ltL^{\mathrm{glob}}(\sigma,E)
 \quad\Longleftrightarrow\quad D^\delta\ltL^\delta E^\delta.
\]
In particular, for filters $D,E$ already on the same ordinal $\delta$,
\[
 \operatorname{can}(D)\ltL^{\mathrm{glob}}\operatorname{can}(E)
 \quad\Longleftrightarrow\quad D\ltL^\delta E,
\]
since extending their canonical presentations back to $\delta$
recovers $D,E$.

\begin{proposition}\label{prop:predecessor-bound}
Let $H$ be a filter on $\delta$ and $E$ a filter on $\eta\leq\delta$.
If $H\ltL^\delta E^\delta$, then $\eta\in H$. Consequently, for
canonical pairs,
\begin{equation}\label{eq:predecessor-ordinal}
 (\rho,D)\ltL^{\mathrm{glob}}(\sigma,E)
 \quad\Longrightarrow\quad \rho\leq\sigma.
\end{equation}
\end{proposition}

\begin{proof}
Let $p$ be a positive reduction witnessing $H\ltL^\delta E^\delta$.
For each $\alpha<\eta$, we have $\delta\cap\alpha=\eta\cap\alpha=\alpha$.
Since $p$ is super-Lipschitz, $p(\delta)\cap\eta=p(\eta)\cap\eta$.
Since $\eta \in E$ and $p$ is a positive reduction, we have
\[
 \begin{aligned}
  \delta\in H &\iff p(\delta)\in E^{\delta} \iff p(\delta)\cap\eta\in E \\
              &\iff p(\eta)\cap \eta \in E \iff p(\eta) \in E \iff \eta\in H.
 \end{aligned}
\]
Since $\delta\in H$, we get $\eta\in H$.

For the second part, let  $(\rho,D)\ltL^{\mathrm{glob}}(\sigma,E)$ and $\tau = \max \{\rho,\sigma\}$. We have $D^{\tau} \ltL^{\tau} E^{\tau}$. Take $H=D^\tau$ and
$\eta=\sigma$. Hence by the first part we have $\sigma\in D^\tau$. So $\rho=\rho(D^\tau)\leq\sigma$,
giving \eqref{eq:predecessor-ordinal}.
\end{proof}

The following adapts the globalization argument of
	\cite[Lemma~3.3.18]{GoldbergBook}, with
	Theorem~\ref{thm:main} supplying local well-foundedness for the
	Lipschitz order.

\begin{corollary}[$\ZF+\DC$]\label{cor:class}
On the class of proper countably complete filters in their canonical
presentations, the global strict Lipschitz order
$\ltL^{\mathrm{glob}}$ is set-like and well-founded. The same holds on
the subclass of countably complete ultrafilters.
\end{corollary}

\begin{proof}
For canonical pairs $(\rho,D)$, $(\sigma,E)$, and $(\xi,H)$, put
$\delta=\max\{\rho,\sigma,\xi\}$. By Proposition~\ref{prop:order-extension},
the global comparisons of these pairs are the comparisons of
$D^\delta,E^\delta,H^\delta$ on $\delta$. Transitivity on $\delta$ therefore
gives transitivity of $\ltL^{\mathrm{glob}}$. Irreflexivity follows from
that of $\ltL^\rho$ for each ordinal $\rho$.

Let $\mathcal C$ be the class of canonical pairs whose filters are
countably complete. For each ordinal $\eta$, put
\[
 \mathcal C_\eta=\{(\rho,D)\in\mathcal C:\rho\leq\eta\}.
\]
This is a set, since
$\mathcal C_\eta\subseteq(\eta+1)\times\Pow(\Pow(\eta))$.
By Proposition~\ref{prop:predecessor-bound}, $\mathcal C_\eta$ contains
every predecessor in $\mathcal C$ of each of its elements.
In particular, all predecessors in $\mathcal C$ of $(\sigma,E)$ lie in
$\mathcal C_\sigma$. Hence the global order is \emph{set-like},
meaning that the predecessors of each element form a set.

The extension map from $\mathcal C_\eta$ to the set of proper
countably complete filters on $\eta$, defined by
\[
 \iota_\eta(\rho,D)=D^\eta,
\]
is a bijection, with inverse $H\mapsto\operatorname{can}(H)$, by
Proposition~\ref{prop:extension-restriction}. By
Proposition~\ref{prop:order-extension}, for $a,b\in\mathcal C_\eta$,
\[
 a\ltL^{\mathrm{glob}}b
 \quad\Longleftrightarrow\quad
 \iota_\eta(a)\ltL^\eta\iota_\eta(b).
\]
Theorem~\ref{thm:main} therefore implies that
$\ltL^{\mathrm{glob}}$ is well-founded on $\mathcal C_\eta$.
To obtain global well-foundedness, let $\mathcal A$ be a nonempty
subclass of $\mathcal C$ and choose $(\sigma,E)\in\mathcal A$.
The nonempty set $\mathcal A\cap\mathcal C_\sigma$ has a minimal
element. That element is also minimal in $\mathcal A$, since
$\mathcal C_\sigma$ is downward closed.

The restriction to countably complete ultrafilters inherits
set-likeness and well-foundedness.
\end{proof}

\begin{proposition}\label{prop:ultrafilters-distinct-ordinals}
Let $(\rho,U)$ and $(\sigma,W)$ be canonical presentations of
ultrafilters. If $\rho<\sigma$, then
$(\rho,U)\ltL^{\mathrm{glob}}(\sigma,W)$.
\end{proposition}

\begin{proof}
Since $\rho\notin W$, we have
$T=\sigma\setminus\rho\in W$. Define $p:\Pow(\sigma)\to\Pow(\sigma)$ by
\[
 p(X)=
 \begin{cases}
  T,&\text{if }X\cap\rho\in U,\\
  \varnothing,&\text{if }X\cap\rho\notin U.
 \end{cases}
\]
For $\alpha<\rho$, every $p(X)$ omits $\alpha$ since $p(X)\subseteq T$.
If $\rho\leq\alpha<\sigma$ and $X\cap\alpha=Y\cap\alpha$, then
$X\cap\rho=Y\cap\rho$ and hence $p(X)=p(Y)$.
Thus $p$ is super-Lipschitz.
Moreover, $p(X)\in W$ if and only if $X\cap\rho\in U$, equivalently
$X\in U^\sigma$. The map $q(X)=\sigma\setminus p(X)$ is a
super-Lipschitz negative reduction, since $W$ is an ultrafilter.
Therefore $U^\sigma\ltL^\sigma W$, or equivalently
$(\rho,U)\ltL^{\mathrm{glob}}(\sigma,W)$.
\end{proof}

For countably complete ultrafilters, this comparison also follows
from the corresponding fact for the Ketonen order
\cite[Lemma~3.3.15]{GoldbergBook}.

\section{Countably incomplete ultrafilters}\label{sec:incomplete-ultrafilters}

Countable completeness plays an essential role in the well-foundedness
result. We now examine descending chains involving countably incomplete
ultrafilters:  in the proof of what follows we combine a regressive pushforward with the shift
	construction of \cite[proof of Theorem~3.4.30]{GoldbergBook},
	adapted to ultrafilters.

\begin{proposition}\label{prop:incomplete}
Every countably incomplete ultrafilter on an ordinal $\delta$ starts
a strictly descending sequence of countably incomplete ultrafilters
on $\delta$.
\end{proposition}

\begin{proof}
Let $W$ be a countably incomplete ultrafilter on $\delta$. Take a
sequence $\{X_n : n<\omega\} \subseteq W$ such that $\bigcap_{n<\omega}X_n\not\in W$. Write $B = \bigcap_{n<\omega}X_n$.
Since $\delta\setminus B\in W$, the sets
$E_n=(\bigcap_{k\leq n}X_k)\setminus B$ form a decreasing sequence in
$W$ with empty intersection. Define $g:\delta\to\omega$ by
$g(\xi)=\min\{n:\xi\notin E_n\}$. Then $\{\xi:g(\xi)>n\}=E_n$.

Since $W$ is nonprincipal, $\delta\geq\omega$ and $W$ contains the
complement of every finite set. Define $f:\delta\to\omega$ by
$f(0)=0$, $f(k)=\min(g(k),k-1)$ for $0<k<\omega$, and $f(\xi)=g(\xi)$
for $\omega\leq\xi<\delta$. Then $f(\xi)<\xi$ for $\xi>0$. For each
$n$, we have  \[f^{-1}[\{ k \in \omega : k >n \}] = \{\xi:f(\xi)>n\} \supseteq E_n\setminus(n+2),\] so $f^{-1}[\{ k \in \omega : k >n \}] \in W$. Thus $D_0=\{X\subseteq\omega:f^{-1}[X]\in W\}$
is a nonprincipal ultrafilter on $\omega$. Define $h:\Pow(\delta)\to\Pow(\delta)$ by
\[
 h(X)=\{\xi<\delta:0<\xi \wedge f(\xi)\in X\}.
\]
For $0<\xi<\delta$, we have $\xi\in h(X)$ if and only if
$f(\xi)\in X\cap\xi$, since $f(\xi)<\xi$. Also, $0$ never belongs to
$h(X)$. Thus $h$ is super-Lipschitz and for every $X \subseteq \delta$ we have  $h(X) = f^{-1}[X\cap\omega]\setminus \{0\}$. So \[h(X)\in W \iff  f^{-1}[X \cap \omega] \in W \iff X \in D_0^\delta.\] Hence $D_0^\delta\ltL^\delta W$.

On $\omega$ let $p(0)=0$ and $p(k)=k-1$ for $k>0$. Set \[D_{n+1}=\{ X \subseteq \omega : p^{-1}[X] \in D_n \}.\] Every $D_n$ is a nonprincipal ultrafilter on $\omega$. The map $X\mapsto\{k+1:k\in X\}$ is super-Lipschitz and for every $X \subseteq \omega$, $\{k+1:k\in X\} = p^{-1}[X] \setminus \{0\}$. Therefore \[ \{k+1:k\in X\} \in D_n \iff p^{-1}[X]\in D_n \iff X \in D_{n+1}\] and
$D_{n+1}\ltL D_n$. By Proposition~\ref{prop:order-extension}, extending
these comparisons to $\delta$ yields the desired sequence starting at
$W$. The extensions of the $D_n$ are
countably incomplete, since the complements in $\delta$ of finite
initial segments of $\omega$ have intersection $\delta\setminus\omega$,
which is omitted by each $D_n^\delta$.
\end{proof}

\section*{Declaration of AI use}

OpenAI Codex was used to explore and develop the mathematical arguments,
check references, review proofs, and draft and revise the manuscript.
The authors take responsibility for the mathematical content and
presentation.


\begin{thebibliography}{9}

\bibitem{Goldberg}
Gabriel Goldberg, \emph{The Ketonen order}, J. Symbolic Logic \textbf{85}
(2020), no.~2, 585--604.
\doi{10.1017/jsl.2020.5}.

\bibitem{GoldbergBook}
Gabriel Goldberg, \emph{The Ultrapower Axiom}, De Gruyter Series in Logic
and Its Applications, vol.~10, De Gruyter, Berlin/Boston, 2022.
\doi{10.1515/9783110719734}.

\bibitem{Kaplan}
Eyal Kaplan, \emph{A note on the Ketonen order and Lipschitz reducibility
	between ultrafilters}, Ann. Pure Appl. Logic \textbf{177} (2026), no.~10,
103805. \doi{10.1016/j.apal.2026.103805}.
Preprint: \href{https://arxiv.org/abs/2512.12835v1}{arXiv:2512.12835v1}.

\bibitem{Ketonen}
Jussi Ketonen, \emph{Strong compactness and other cardinal sins},
Ann. Math. Logic \textbf{5} (1972), no.~1, 47--76.
\doi{10.1016/0003-4843(72)90018-6}.

\bibitem{Mitchell}
William J. Mitchell, \emph{Sets constructible from sequences of
	ultrafilters}, J. Symbolic Logic \textbf{39} (1974), no.~1, 57--66.
\doi{10.2307/2272343}.

\bibitem{Wadge}
William W. Wadge, \emph{Reducibility and determinateness on the Baire
	space}, Ph.D. thesis, University of California, Berkeley, 1983.
\url{https://billwadge.com/wp-content/uploads/2021/02/wadgephd.pdf}.

\bibitem{Yablo}
Stephen Yablo, \emph{Paradox without self-reference}, Analysis \textbf{53}
(1993), no.~4, 251--252.
\doi{10.1093/analys/53.4.251}.

\end{thebibliography}
\end{document}